 \documentclass[12pt,reqno]{article}

 \usepackage{amsmath, amsthm, amsfonts}
 \usepackage{bm}
 \usepackage{mathrsfs, amssymb}
 \usepackage{graphicx, color, bm}

 \numberwithin{equation}{section}

 \newtheorem{thm0}{Theorem}[section]
 \newtheorem{exa0}{Theorem}[section]

 \newtheorem{def1}[thm0]{Definition}
 \newtheorem{lem1}[thm0]{Lemma}
 \newtheorem{thm1}[thm0]{Theorem}
 \newtheorem{cor1}[thm0]{Corollary}
 \newtheorem{pro1}[thm0]{Proposition}
 \newtheorem{con1}[thm0]{Condition}
 \newtheorem{rem1}[thm0]{Remark}
 \newtheorem{exa1}[exa0]{\it{Example}}

 \def\bglemma{\begin{lem1}}\def\edlemma{\end{lem1}}
 \def\bgtheorem{\begin{thm1}}\def\edtheorem{\end{thm1}}
 
 \def\bgproposition{\begin{pro1}}\def\edproposition{\end{pro1}}

 \def\benumerate{\begin{enumerate}}\def\eenumerate{\end{enumerate}}
 \def\bitemize{\begin{itemize}}\def\eitemize{\end{itemize}}

 \def\beqlb{\begin{eqnarray}}\def\eeqlb{\end{eqnarray}}
 \def\beqnn{\begin{eqnarray*}}\def\eeqnn{\end{eqnarray*}}
\def\E{\mathbf{E}}\def\p{\mathbf{P}}\def\d{{\mbox{\rm d}}}\def\e{{\mbox{\rm e}}}
 \def\mcr{\mathscr}\def\mbb{\mathbb}\def\z{{\bm z}}
 \def\X{{\bm X}}\def\x{{\bm x}}\def\y{{\bm y}}\def\b{{\bm b}}
 \def\D{{\mbb D}} \def\mrm{\mathrm}

 \def\eqref#1{{\rm(\ref{#1})}}

\def\ar{\!\!\!&}

\begin{document}


\thispagestyle{empty}
\bigskip
\centerline{\bf\LARGE Moment properties for two-type }

\medskip
\centerline{\bf\LARGE  continuous-state branching processes in}

\medskip
\centerline{\bf\LARGE   L\'evy random environments}

\bigskip

\centerline{Shukai Chen and Xiangqi Zheng*}

\medskip

\centerline{\it School of Mathematics and Statistics,Fujian Normal University}

\centerline{\it Fuzhou, 350007, People's Republic of China}

\centerline{\it School of Mathematics, East China University of Science and Technology,}

\centerline{\it Shanghai, 200237, People's Republic of China}

\smallskip

\centerline{\it E-mail: skchen@mail.bnu.edu.cn and zhengxq@ecust.edu.cn}

\bigskip

\bigskip

\noindent{\bf Abstract.} We first derive the recursions for integer moments of two-type continuous-state branching processes in L\'{e}vy random environments. We show that the $n$th moment of the process is a polynomial of the initial value of the process with at most $n$ degree. Meanwhile, the criteria for the existence of $f$-moment of the process is also established under some natural conditions.

\bigskip


\bigskip

\section{Introduction and main results}

The branching model in continuous time and state called {\it continuous
state branching process} (CB-process) was first introduced by  Ji\v{r}ina \cite{Ji58}. This model can be regarded as the scaling limits of classical Galton-Watson branching processes, see \cite{ref2,ref3,ref4}. Suppose that $(\Omega,{\cal{F}},{\cal{F}}_t,{\bf P})$ is a filtered probability space satisfying the usual hypotheses. Given $\sigma\geq0, b\in\mbb{R}$ constants. Let $m(\d z)$ be a $\sigma$-finite measure on $[0,\infty)$ satisfying $\int_0^{\infty}(z\wedge z^2)\,m(\d z)<\infty.$ Suppose that $\{B(t)\}$ is a standard $({\cal{F}}_t)$-Brownian motion and $\{M(\d s,\d z,\d u)\}$ is a $({\cal{F}}_t)$-Poisson random measure with intensity $\d sm(\d z)\d u$ and $\tilde{M}(\d s,\d z,\d u)=M(\d s,\d z,\d u)-\d sm(\d z)\d u$ is the compensated measure. We assume those two noises are independent. By Theorem 3.1 in \cite{ref6}, there is a unique
positive strong solution to
\beqlb\label{cb}
&&X(t)= X(0)+\int_0^t\sqrt{2\sigma X(s)}\d B(s)+\int_0^t\int_0^{1}\int_0^{X(s-)}z\tilde{M}(\d s,\d z,\d u)\nonumber\\
&&~~~~~~~~~~~~~~~~
-\int_0^tbX(s)\d s+\int_0^t\int_1^{\infty}\int_0^{X(s-)}z\,M(\d s,\d z,\d u)
\eeqlb
and the solution $\{X(t):t\geq0\}$ is a CB-process with branching mechanism $\phi$ satisfying a L\'evy-Khintchine's type
$$
\phi(\lambda)=b\lambda+\sigma\lambda^2+\int_0^{\infty}(\e^{-\lambda z}-1+\lambda z{\bf1}_{\{z\leq1\}})\,m(\d z),\quad\quad \lambda\geq0.
$$
Based on \eqref{cb}, \cite{ref6} studied the moment property of CB-processes. Indeed, suppose that $f$ is a positive continuous function on
$[0, \infty)$ satisfying

\smallskip

\noindent\textbf{Condition A.~} There exist constants $c\ge 0$ and $K> 0$ such that

(A1) $f$ is convex on $[c,\infty)$;

(A2) $f(xy)\le Kf(x)f(y)$ for all $x,y\in [c,\infty)$.

\noindent If $\{X(t):t\ge 0\}$ is the strong solution of equation \eqref{cb} with
$\p(X(0) > 0) > 0,$ then for any $t > 0,$ the equivalent condition of $\E f(X(t)) < \infty$ is $\E f(X(0)) < \infty$ and $ \int_1^{\infty}f(z)m(dz) < \infty.$ Later, \cite{ref7} generalised the result to {\it continuous-state branching processes in L\'{e}vy random environments} (CBRE-processes). A CBRE-process $\{Y(t):t\geq0\}$ can be seen as the unique strong solution of
\begin{eqnarray}\label{cbre}
Y(t)\ar=\ar Y(0)-\int_0^tbY(s)\d s+\int_0^t\sqrt{2\sigma Y(s)}\d B(s)+
\int_0^t\int_1^{\infty}\int_0^{Y(s-)}zM(\d s,\d z,\d u)\cr\cr
\ar\ar\quad\quad+\int_0^t\int_0^{1}\int_0^{Y(s-)}z\tilde{M}(\d s,\d z,\d u)+\int_0^t Y(s-)\d L(s),
\end{eqnarray}
where $b,\sigma,B,M$ are the same as that in \eqref{cb}, $\{L(t):t\ge 0\}$ is a $(\mathscr{F}_t)$-L\'{e}vy process defined by
 \begin{eqnarray}\label{L}
 L(t) = \beta t + \sigma_1B(t) + \int_0^t \int_{\D_1}(e^z - 1)\tilde{N}(\d s, \d z)
+ \int_0^t \int_{\D_1^c} (e^{z} - 1) N(\d s, \d z),
 \end{eqnarray}
where $\D_1=[-1, 1],$ $\beta \in \mathbb{R},\sigma_1\geq 0$, $\{B(t)\}$~is a standard $({\cal{F}}_t)$-Brownian motion and $\{N(\d s, \d z)\}$ is a $({\cal{F}}_t)$-Poisson random measure on $(0, \infty)\times \mathbb{R}$ with intensity $\d s \nu(\d z),$ $\nu$ is a $\sigma$-finite measure satisfying $\int_0^{\infty}(1\wedge z^2) \nu(\d z)<\infty.$ In particular, when $L(t)\equiv0$ for all $t\geq0$, the CBRE-process reduces to a CB-process with branching mechanism $\phi$. An associated L\'evy process $\{\xi(t):t\geq0\}$ is defined by
\beqlb\label{xit}
\xi(t)=at+\sigma_1B(t)+\int_0^t \int_{\D_1}z\tilde{N}(\d s, \d z)
+ \int_0^t \int_{\D_1^c} z N(\d s, \d z),
\eeqlb
where $a=\beta-\frac{1}{2}{\sigma_1}^2-\int_{\D_1}(\e^z-1-z)\,\nu(\d z)$. Clearly, the two processes $\{\xi(t):t\geq0\}$ and $\{L(t):t\geq0\}$ generate the same
filtration. We refer to \cite{ref8,ref9} for more details of CBRE-processes. Under the hypothesis $\p(Y(0)> 0)> 0,$ \cite{ref7} showed that $\E f(Y(t))< \infty$ for any $t>0$ if and only if $\E f(Y(0))< \infty$, $\int_1^\infty f(z) m(\d z)< \infty$ and $\int_1^\infty f(e^z) \nu(\d z) < \infty,$ here the function $f$ satisfies {\bf Condition A}. A recursion formula of the $n$-moment of such processes is also established in \cite{ref7}.

A {\it two-type continuous-state branching process in L\'evy random environment} (two-type CBRE-process) with a slightly stronger moment condition on the branching mechanism was constructed by \cite{qinzheng}, where the environment process is still defined by \eqref{L} or \eqref{xit}, and the branching mechanism\ $\bm{\phi}=(\phi_1,\phi_2)$\ is a function from\
$\mathbb{R}_+^2$\ to itself with the following representations,
\beqnn
&&\phi_1(\bm{\lambda}) = b_{11}\lambda_1 + b_{12}\lambda_2 + c_1 \lambda_1^2 +
\int_{\mathbb{R}_+^2} (\e^{-\langle\bm{\lambda,z}\rangle} - 1 + \lambda_1z_1)m_1(\d
\z),\\
&&\phi_2(\bm{\lambda}) = b_{21}\lambda_1 + b_{22}\lambda_2 + c_2 \lambda_2^2 +
\int_{\mathbb{R}_+^2} (\e^{-\langle\bm{\lambda,z}\rangle} - 1 + \lambda_2z_2)m_2(\d
\z).
 \eeqnn
Here,\ $(b_{ij})$~is a~$(2\times2)$-matrix with~$b_{12},b_{21}\le 0,$ $c_1,c_2\ge0$, $m_1,m_2$~are~$\sigma$-finite measures on~$\mathbb{R}_+^2$~supported
by~$\mathbb{R}_+^2\setminus\{\mathbf{0}\},$ satisfying
 \begin{eqnarray*}
\int_{\mathbb{R}_+^2\setminus\{\bf0\}} (z_1\wedge z_1^2 + z_2)m_1(\d\z) + \int_{\mathbb{R}_+^2\setminus\{\bf0\}}
(z_2\wedge z_2^2 + z_1)m_2(\d\z) < \infty.
 \end{eqnarray*}
A two-dimensional Markov process $\{\X(t)=(X_1(t),X_2(t)):t\ge 0\}$ is a CBRE-process if its transition semigroup $\{P_t:t\geq0\}$ is determined by
\beqlb\label{Pt}
\int_{\mbb{R}_+^2}\mrm{e}^{-\langle{\bm\lambda},{\bm y}\rangle}\,P_t({\bm x},\d{\bm y})=
\E\Big[\mrm{exp}\Big\{-\langle{\bm x},{\bm v}_{0,t}(\xi,{\bm\lambda})\rangle\Big\}\Big],
\eeqlb
where $r\mapsto{\bm v}_{r,t}(\xi,{\bm \lambda})=(v^{(1)}_{r,t}(\xi,{\bm\lambda}),v^{(2)}_{r,t}(\xi,{\bm\lambda}))$ is the unique solution to
$$
v^{(i)}_{r,t}(\xi,{\bm\lambda})=\e^{\xi(t)-\xi(r)}\lambda_i
-\int_r^t\e^{\xi(s)-\xi(r)}\phi_i({\bm v}_{s,t}(\xi,{\bm \lambda}))\,\d s,\quad i=1,2,\quad r\leq t\in\mbb{R}
$$
and we take $r=0$ in \eqref{Pt}. 

Let $\{B_1(t)\}$ and $\{B_2(t)\}$ be two standard $({\cal{F}}_t)$-Brownian motions, $\{M_1(\d s,\d u,\d\z)\}$ and $\{M_2(\d s,\d u,\d\z)\}$ be two $({\cal{F}}_t)$-Poisson random measures on $(0,\infty)^4$ with characteristic measures $m_1$ and $m_2,$ respectively, $\{\tilde{M}_1(\d s,\d u,\d\z)\}$ and $\{\tilde{M}_2(\d s,\d u,\d\z)\}$ are associated compensated measures. We assume those random elements are mutually independent. It was proved in \cite{qinzheng} that a two-type CBRE-process with branching mechanism\ $\bm{\phi}=(\phi_1,\phi_2)$\ and environment $\{\xi(t):t\geq0\}$ or $\{L(t):t\geq0\}$ can also be seen as the unique non-negative strong solution to
  \beqlb
&&X_1(t)=X_1(0)-\int_0^t\Big(b_{11}X_1(s)+b_{21}X_2(s)\Big)\d s
+ \int_0^t\sqrt{2c_1X_1(s)}\d B_1(s)\nonumber\\
&&~~~~~~~
+\int_0^t\int_0^{X_1(s-)}\int_{\mathbb{R}_+^2\setminus\{\mathbf{0}\}}z_1\tilde{M}_1(\d
s,\d u,\d\z)+\int_0^tX_1(s-)\d L(s)\nonumber\\
&&~~~~~~~
+\int_0^t\int_0^{X_2(s-)}\int_{\mathbb{R}_+^2\setminus\{\mathbf{0}\}}z_1M_2(\d
s,\d u,\d\z),\label{eq01}\\
&&X_2(t)=X_2(0)-\int_0^t\Big(b_{12}X_1(s)+b_{22}X_2(s)\Big)\d s+
 \int_0^t\sqrt{2c_2X_2(s)}\d B_2(s)\cr\cr
&&~~~~~~~
+\int_0^t\int_0^{X_1(s-)}\int_{\mathbb{R}_+^2\setminus\{\mathbf{0}\}}z_2M_1(\d
s,\d u,\d\z)+\int_0^tX_2(s-)\d L(s)\cr\cr
&&~~~~~~~
+\int_0^t\int_0^{X_2(s-)}\int_{\mathbb{R}_+^2\setminus\{\mathbf{0}\}}z_2\tilde{M}_2(\d
s,\d u,\d\z).\label{eq02}
\eeqlb
In particular, when $L(t)\equiv0$ for all $t\geq0$, the two-type CBRE-process reduces to a two-type CB-process with branching mechanism $\bm{\phi}$. One can refer to \cite{bar},\cite{machunhua},\cite{marugang} for more details on the multi-type continuous-state branching processes. Moreover, using the tool of stochastic differential equations, Barczy, Li and Pap\cite{barm} calculate the $n$-moment of multi-type continuous-state branching process. Our result is an extension of \cite{barm}. Let~$\{\X(t)=(X_1(t),X_2(t)):t\ge 0\}$~be the strong solution of the stochastic equation system \eqref{eq01}--\eqref{eq02}. We calculate the integer moment recursions of $\{\X(t)=(X_1(t),X_2(t)):t\ge 0\}$ and the equivalent condition for the existence of $f$-moment. 

\bgtheorem
{\bf(n-moment)}

{\it Suppose that there exists an integer $n\geq2$ such that
$$
\E \|\X(0)\|^n<\infty,~~\int_{\mathbb{R}_+^2\setminus\{{\bf 0}\}}\|\z\|^n(m_1+m_2)(\d \z)<\infty, ~~\int_1^{\infty}\e^{nz}\nu(\d z)<\infty.
$$
Then,
\begin{eqnarray*}
 \E[X_1(t)]^n\ar=\ar\E [X_1(0)]^n\E\e^{n\bar{\xi}_1(t)} + \sum_{j = 0}^{n - 2}A^1_{n,j} \int_0^t\E X_1(s-)^{j+1} \E\e^{n\bar{\xi}_1(t)-n\bar{\xi}_1(s-)}\d s\cr\cr
\ar\quad\ar+ \sum_{j = 0}^{n - 1}B^1_{n,j} \int_0^t\E X_1(s-)^j X_2(s-)\E\e^{n\bar{\xi}_1(t)-n\bar{\xi}_1(s-)}\d s
\end{eqnarray*}
and}
\begin{eqnarray*}
 \E[X_2(t)]^n\ar=\ar\E [X_2(0)]^n\E\e^{n\bar{\xi}_2(t)} + \sum_{j = 0}^{n - 2}A^2_{n,j} \int_0^t\E X_2(s-)^{j+1} \E\e^{n\bar{\xi}_2(t)-n\bar{\xi}_2(s-)}\d s\cr\cr
\ar\quad\ar+ \sum_{j = 0}^{n - 1}B^2_{n,j} \int_0^t\E X_2(s-)^j X_1(s-)\E\e^{n\bar{\xi}_2(t)-n\bar{\xi}_2(s-)}\d s,
\end{eqnarray*}
{\it where $\bar{\xi}_1(t)=\xi(t)-b_{11}t,\,\bar{\xi}_2(t)=\xi(t)-b_{22}t$},
\begin{eqnarray*}
A^1_{n,j}\ar=\ar \binom{n}{j}\int_{\mathbb{R}_+^2\setminus\{\mathbf{0}\}} z_1^{n - j}m_1(\d
\z),\qquad 0\le j<n-2,\cr\cr
A^1_{n,j}\ar=\ar \binom{n}{n-2}\int_{\mathbb{R}_+^2\setminus\{\mathbf{0}\}} z_1^{n - j}m_1(\d \z)+c_1n(n-1),\qquad j=n-2,\cr\cr
B^1_{n,j}\ar=\ar \binom{n}{j}\int_{\mathbb{R}_+^2\setminus\{\mathbf{0}\}} z_1^{n - j}m_2(\d \z),\qquad j<n-1,\cr\cr
B^1_{n,j}\ar=\ar \binom{n}{n-1}\int_{\mathbb{R}_+^2\setminus\{\mathbf{0}\}} z_1^{n - j}m_2(\d \z)-b_{21}n,\qquad j=n-1,\cr\cr
A^2_{n,j}\ar=\ar \binom{n}{j}\int_{\mathbb{R}_+^2\setminus\{\mathbf{0}\}} z_2^{n - j}m_2(\d \z),\qquad j<n-2,\cr\cr
A^2_{n,j}\ar=\ar \binom{n}{n-2}\int_{\mathbb{R}_+^2\setminus\{\mathbf{0}\}} z_2^{n - j}m_2(\d \z)+c_2n(n-1),\qquad j=n-2,\cr\cr
B^2_{n,j}\ar=\ar \binom{n}{j}\int_{\mathbb{R}_+^2\setminus\{\mathbf{0}\}} z_2^{n - j}m_1(\d \z),\qquad j<n-1,\cr\cr
B^2_{n,j}\ar=\ar \binom{n}{n-1}\int_{\mathbb{R}_+^2\setminus\{\mathbf{0}\}} z_2^{n - j}m_1(\d \z)-b_{12}n,\qquad j=n-1.
\end{eqnarray*}
{\it Moreover, for each $t\ge 0,k=1,2,\cdots,n$ and $i=1,2,$ there exists a polynomial function $Q_{i,t,k}:\mathbb{R}_+^2\rightarrow\mathbb{R}_+$ having degree at most $k$ such that,
$$\E [X_i(t)]^k=\E[Q_{i,t,k}(\X(\mathbf{0}))].$$}
\edtheorem

\medskip

\bgtheorem
{\bf(}$\bm f$-{\bf moment)}\label{thm2020021401}

{\it Suppose that $\p (\|\X(0)\|>0)>0.$ Suppose further that $f$ satisfies \textbf{Condition A}. Then for any $t>0,$ $\E f(\|\X(t)\|)<\infty$ if and only if}
$$
\E f(\|\X(0)\|)<\infty,~\int_{\{\|\z\|\ge 1\}}f(\|\z\|)(m_1+m_2)(\d\z)<\infty, ~\int_1^{\infty}f(\e^z)\nu(\d z)<\infty.
$$
\edtheorem
 The proofs will be given in the next two sections.

{\bf Notation}. For any $\{\x(t):t\geq0\}$ taking values in $\mathcal{S}$, the space of two-dimensional functions with c\'adl\'ag paths, write $\x(t)=(x_1(t),x_2(t))$. For any $\x=(x_1,x_2),\y=(y_1,y_2)\in\mbb{R}_{+}^2$, write $\y\geq\x$ if $y_1\geq x_1$ and $y_2\geq x_2$. Throughout this paper, we make the conventions
$$
\int_a^b=\int_{(a,b]}\,\,\,\,\mbox{and}\,\,\,\,\int_a^{\infty}=\int_{(a,\infty)}
$$
for any $b\geq a\geq0$. Given a function $f$ defined on a subset of $\mbb{R}$, we write
$$
\Delta_zf(x)=f(x+z)-f(x)\,\,\,\,\mbox{and}\,\,\,\,D_zf(x)=\Delta_zf(x)-f'(x)z
$$
for $x,z\in\mbb{R}$ if the right-hand side is meaningful.

\section{n-moment of two-type CBRE processes}

\setcounter{equation}{0}

In this section, we derive the recursions of integer moments of the process with the help of truncated processes $\{\X^{(k)}(t):t\ge 0\}$$(k=1,2,\cdots)$ defined as the unique non-negative strong solution of the following equation system

  \begin{eqnarray*}
&&X^{(k)}_1(t)=X^{(k)}_1(0) + \int_0^t\sqrt{2c_1X^{(k)}_1(s)}\d B_1(s) -
\int_0^t\Big(b_{11}X^{(k)}_1(s)+b_{21}X^{(k)}_2(s)\Big)\d s
\cr\cr
&&~~~~~~~~~
+\int_0^t\int_0^{X^{(k)}_1(s-)}\int_{\mathbb{R}_+^2\setminus\{\mathbf{0}\}}z_1\mathbf{1}_{\{\|\z\|
\leq k\}}\tilde{M}_1(\d s,\d u,\d\z)+\int_0^tX^{(k)}_1(s-)\d L^{(k)}(s)\cr\cr
&&~~~~~~~~~
+\int_0^t\int_0^{X^{(k)}_2(s-)}\int_{\mathbb{R}_+^2\setminus\{\mathbf{0}\}}z_1\mathbf{1}_{\{\|\z\|
\leq k\}}M_2(\d s,\d u,\d\z),
\end{eqnarray*}
  \begin{eqnarray*}
&&X^{(k)}_2(t)=X^{(k)}_2(0) + \int_0^t\sqrt{2c_2X^{(k)}_2(s)}\d B_2(s) -
\int_0^t\Big(b_{12}X^{(k)}_1(s)+b_{22}X^{(k)}_2(s)\Big)\d s\cr\cr
&&~~~~~~~~~
+\int_0^t\int_0^{X^{(k)}_1(s-)}\int_{\mathbb{R}_+^2\setminus\{\mathbf{0}\}}z_2\mathbf{1}_{\{\|\z\|
\leq k\}}M_1(\d s,\d u,\d\z)+\int_0^tX^{(k)}_2(s-)\d L^{(k)}(s)\cr\cr
&&~~~~~~~~~
+\int_0^t\int_0^{X^{(k)}_2(s-)}\int_{\mathbb{R}_+^2\setminus\{\mathbf{0}\}}z_2\mathbf{1}_{\{\|\z\|
\leq k\}}\tilde{M}_2(\d s,\d u,\d\z),
\end{eqnarray*}
where
$\{L^{(k)}(t): t \geq 0\}$ is a L\'{e}vy processes
with the following L\'{e}vy-It\^{o} decomposition,
$$L^{(k)}(t) = \beta t + \sigma_1 W(t) + \int_0^t \int_{\D_1}(e^z -
1)\tilde{N}(\d s, \d z)+ \int_0^t \int_{\D_1^c} \left(e^{z1_{\left\{z \leq
k\right\}}} - 1\right) N(\d s, \d z).
$$
The associated L\'evy process $\{\xi^{(k)}(t):t\geq0\}$ is defined by
$$
\xi^{(k)}(t) = at + \sigma_1 W(t) + \int_0^t \int_{\D_1}z\tilde{N}(\d s, \d z) + \int_0^t \int_{\D_1^c} z1_{\left\{z \leq
k\right\}} N(\d s, \d z).
$$
As $k\rightarrow\infty$, the truncated sequence $\{\X^{(k)}(t):t\ge 0\}_{k\geq1}$ converges increasingly to $\{\X(t):t\ge 0\}.$ And we will prove it through the following two propositions.
 \bgproposition
{\it For all $n_2\ge n_1>1,$}
$$
\p\Big(\X^{(n_1)}(t) \leq \X^{(n_2)}(t)\,\, \text{for\ all}\ t \geq 0\Big) = 1.
$$
\edproposition

 \proof
Fix $n_1,n_2$ with $n_2\ge n_1>1$. Define~${\bm{\zeta}}(t)=\X^{(n_1)}(t)-\X^{(n_2)}(t).$
\beqnn
\zeta_1(t) \ar\ \le\ \ar \zeta_1(0) +
\int_0^t(\sqrt{2c_1X^{(n_1)}_1(s)}-\sqrt{2c_1X^{(n_2)}_1(s)})\d B_1(s)\cr\cr
\ar\ \ar + \int_0^t[-b_{11}\zeta_1(s)-b_{21}\zeta_2(s)]\d s
\cr\cr
&
&+\int_0^t\int_{X^{(n_2)}_1(s-)}^{X^{(n_1)}_1(s-)}\int_{\mathbb{R}_+^2\setminus\{\mathbf{0}\}}z_1\mathbf{1}_{\{\|\z\|
\leq {n_2}\}}\tilde{M}_1(\d s,\d u,\d\z)\cr\cr
&
&+\int_0^t\int_{X^{(n_2)}_2(s-)}^{X^{(n_1)}_2(s-)}\int_{\mathbb{R}_+^2\setminus\{\mathbf{0}\}}z_1\mathbf{1}_{\{\|\z\|
\leq n_2\}}M_2(\d s,\d u,\d\z)\cr\cr
& &+\int_0^t\zeta_1(s-)\d L^{(n_2)}(s).
\eeqnn
For $m\geq1$, define~$\tau_m:=\inf\{t\ge 0:X^{(n_1)}_1(t)\vee X^{(n_2)}_1(t)\vee
X^{(n_1)}_2(t)\vee X^{(n_2)}_2(t)\ge m\}.$ For~$t\ge 0,$~choose a
decreasing sequence~$\{a_k\}$~such that $a_0=1,$~$a_k\rightarrow
0,$~$\int_{a_k}^{a_{k-1}}z^{-1}\d z=k$ for $k\ge 1.$ Let~$x\mapsto\psi_k(x)$~be
non-negative functions on~$\mathbb{R}$~supported by~$(a_k,a_{k-1}),$~and
satisfying~$\int_{a_k}^{a_{k-1}}\psi_k(x)\d
x=1,$~$0\le\psi_k(x)\le\frac{2}{kx},$~$a_k\le x\le a_{k-1}.$~For~$k\ge
1,$~define a twice-differentiable non-negative function
$$
\varphi_k(z)=\int_0^{z}\d y\int_0^y\psi_k(x)\d x,\quad z\in\mathbb{R}.
$$
By It\^{o}'s formula,
\beqnn
\ar\ar\quad \varphi_k[\zeta_1(t\wedge\tau_m)]\cr\cr
\ar\ar\leq\varphi_k[\zeta_1(0)]
+
\int_0^{t\wedge\tau_m}\varphi'_k[\zeta_1(s)][-b_{11}\zeta_1(s)-b_{21}\zeta_2(s)]\d s
\cr\cr
\ar\ar\quad+
\int_0^{t\wedge\tau_m}\varphi''_k[\zeta_1(s)]c_1(\sqrt{X^{(n_1)}_1(s)}-\sqrt{X^{(n_2)}_1(s)})^2\d
s+\frac{1}{2}\int_0^{t\wedge\tau_m}\varphi_k''(\zeta_1(s))\sigma^2\zeta_1^2(s)\d s\cr\cr
\ar\ar\quad+\int_0^{t\wedge\tau_m}\int_{\mathbb{R}_+^2\setminus\{\mathbf{0}\}}\Delta_{z_1}\varphi_k(\zeta_1(s-))\zeta_2(s-)\mathbf{1}_{\{\|\z\|
\leq n_2\}}\mathbf{1}_{\{\zeta_2(s-)>0\}}m_2(\d\z)\d s\cr\cr
\ar\ar\quad-\int_0^{t\wedge\tau_m}\int_{\mathbb{R}_+^2\setminus\{\mathbf{0}\}}\Delta_{(-z_1)}\varphi_k(\zeta_1(s-))\zeta_2(s-)\mathbf{1}_{\{\|\z\|
\leq n_2\}}\mathbf{1}_{\{\zeta_2(s-)\le 0\}} m_2(\d \z)\d s\cr\cr
\ar\ar\quad+\int_0^{t\wedge\tau_m}\int_{\mathbb{R}_+^2\setminus\{\mathbf{0}\}}D_{z_1}\varphi_k(\zeta_1(s-))\zeta_1(s-)\mathbf{1}_{\{\|\z\|
\leq n_2\}}\mathbf{1}_{\{\zeta_1(s-)>0\}}m_1(\d \z)\d s\cr\cr
\ar\ar\quad-\int_0^{t\wedge\tau_m}\int_{\mathbb{R}_+^2\setminus\{\mathbf{0}\}}D_{(-z_1)}\varphi_k(\zeta_1(s-))\zeta_1(s-)\mathbf{1}_{\{\|\z\|
\leq n_2\}}\mathbf{1}_{\{\zeta_1(s-)\le 0\}}m_1(\d \z)\d s\cr\cr
\ar\ar\quad+\int_0^{t\wedge\tau_m}\int_{\D_1}D_{\zeta_1(s-)(e^z-1)}\varphi_k(\zeta_1(s-))
\nu(\d z)\d s\cr\cr
\ar\ar\quad+ \int_0^{t\wedge\tau_m}\int_{\D_1^c}\Delta_{\zeta_1(s-)(\e^z-1)}\varphi_k(\zeta_1(s-))\mathbf{1}_{\{\|\z\|
\leq n_2\}}\nu(\d z)\d s+mart.
\eeqnn
It is not hard to see that, as\ $k\rightarrow\infty,$
$$\varphi_k''(\zeta_1(s))(\zeta_1(s))^2\rightarrow 0,$$
$$\varphi_k''(\zeta_1(s))\Big(\sqrt{X_1^{(n_1)}(s)}-\sqrt{X^{(n_2)}_1(s)}\Big)^2\rightarrow
0.$$
Moreover,
$$0\le\Delta_{z_1}\varphi_k(\zeta_1(s-))\zeta_2(s-)\mathbf{1}_{\{\zeta_2(s-)>
0\}}\le z_1\zeta_2(s-)^+.$$ Since\ $z\mapsto\varphi(z)$\ is nondecreasing,
$$\Delta_{(-z_1)}\varphi_k(\zeta_1(s-))\zeta_2(s-)\mathbf{1}_{\{\zeta_2(s-)\le
0\}}\ge 0.$$
By Taylor's expansion, when\ $\zeta_1(s-)>0,$
\beqnn
D_{z_1}\varphi_k(\zeta_1(s-))\ar =\ar
z_1^2\int_0^1\varphi_k''(\zeta_1(s-)+tz_1)(1-t)\d t\cr\cr
\ar =\ar z_1^2\int_0^1\psi_k(\zeta_1(s-)+tz_1)(1-t)\d t\cr\cr
\ar \le \ar z_1^2\int_0^1\frac{2}{k(\zeta_1(s-)+tz_1)}(1-t)\d t\cr\cr
\ar \le \ar \frac{z_1^2}{k\zeta_1(s-)}.
\eeqnn
Meanwhile,
$$D_{z_1}\varphi_k(\zeta_1(s-))\le
\varphi_k(\zeta_1(s-)+z_1)-\varphi_k(\zeta_1(s-))\le z_1.$$
Thus,
\beqnn
\ar\ar\int_{\mathbb{R}_+^2\setminus\{\mathbf{0}\}}D_{z_1}\varphi_k(\zeta_1(s-))\zeta_1(s-)\mathbf{1}_{\{\|\z\|
\leq n_2\}}\mathbf{1}_{\{\zeta_1(s-)> 0\}}m_1(\d\z)\cr\cr
\ar=\ar\int_{\mathbb{R}_+^2\setminus\{\mathbf{0}\}}D_{z_1}\varphi_k(\zeta_1(s-))\zeta_1(s-)
\mathbf{1}_{\{\|\z\|\leq n_2,z_1> 1\}}\mathbf{1}_{\{\zeta_1(s-)> 0\}}m_1(\d\z)\cr\cr
\ar\ar\qquad+\int_{\mathbb{R}_+^2\setminus\{\mathbf{0}\}}D_{z_1}\varphi_k(\zeta_1(s-))\zeta_1(s-)\mathbf{1}_{\{\|\z\|
\leq n_2,z_1\le 1\}}\mathbf{1}_{\{\zeta_1(s-)> 0\}}m_1(\d\z)\cr\cr
\ar\le\ar
\int_{\mathbb{R}_+^2\setminus\{\mathbf{0}\}}z_1\zeta_1(s-)\mathbf{1}_{\{\|\z\|
\leq n_2\}}\mathbf{1}_{\{\zeta_1(s-)> 0\}}\mathbf{1}_{\{z_1> 1\}}m_1(\d\z)\cr\cr
\ar\ar\qquad+
\int_{\mathbb{R}_+^2\setminus\{\mathbf{0}\}}\frac{z_1^2}{k\zeta_1(s-)}\zeta_1(s-)\mathbf{1}_{\{\|\z\|
\leq n_2\}}\mathbf{1}_{\{\zeta_1(s-)> 0\}}\mathbf{1}_{\{z_1\le 1\}}m_1(\d\z)\cr\cr
\ar=\ar
\zeta_1(s-)^+\int_{\mathbb{R}_+^2\setminus\{\mathbf{0}\}}
z_1\mathbf{1}_{\{\|\z\|\leq n_2,z_1> 1\}}m_1(\d\z)
+\frac{1}{k}\int_{\mathbb{R}_+^2\setminus\{\mathbf{0}\}}z_1^2\mathbf{1}_{\{\|\z\|
\leq n_2,z_1\le 1\}}m_1(\d\z).
\eeqnn
Similarly, when\ $z\in\D_1,\ \zeta_1(s-)>0,$
\beqnn
\ar\ar D_{\zeta_1(s-)(e^z-1)}\varphi_k(\zeta_1(s-))\cr\cr
\ar\ar\qquad =
\zeta_1^2(s-)(e^z-1)^2\int_0^1\psi_k[\zeta_1(s-)\big(t(e^z-1)+1\big)](1-t)\d
t\cr\cr
\ar\ar\qquad
\le\zeta_1^2(s-)(e^z-1)^2\int_0^1\frac{2(1-t)}{\zeta_1(s-)[t(e^z-1)+1]}\d
t\cr\cr
\ar\ar\qquad \le\zeta_1(s-)(e^z-1)^2.
\eeqnn
Hence,
\beqnn
\ar\ar
\int_0^{t\wedge\tau_m}\int_{\D_1}D_{\zeta_1(s-)(e^z-1)}\varphi_k(\zeta_1(s-))\nu(\d z)\d s\cr\cr
\ar\ar\qquad
\le\int_0^{t\wedge\tau_m}\int_{\D_1}\zeta_1(s-)(e^z-1)^2\nu(\d z)\d s.
\eeqnn
Then,
\beqnn
\varphi_k[\zeta_1(t\wedge\tau_m)]\ar\le\ar
\varphi_k(\zeta_1(0))+\int_0^{t\wedge\tau_m}\varphi'_k[\zeta_1(s)][-b_{11}\zeta_1(s)-b_{21}\zeta_2(s)]\d
s\cr\cr
\ar +\ar
\int_0^{t\wedge\tau_m}\varphi_k''(\zeta_1(s))c_1(\sqrt{X^{(n_1)}_1(s)}-\sqrt{X^{(n_2)}_1(s)})^2\d
s\cr\cr
\ar+\ar\frac{1}{2}\int_0^{t\wedge\tau_m}\varphi_k''(\zeta_1(s))\sigma^2\zeta_1^2(s)\d
s\cr\cr
\ar+\ar\int_0^{t\wedge\tau_m}\int_{\mathbb{R}_+^2\setminus\{\mathbf{0}\}}z_1\zeta_2(s-)^+\mathbf{1}_{\{\|\z\|
\leq n_2\}}m_2(\d \z)\d s\cr\cr
\ar+\ar\int_0^{t\wedge\tau_m}\int_{\mathbb{R}_+^2\setminus\{\mathbf{0}\}}z_1\zeta_1(s-)^+\mathbf{1}_{\{\|\z\|
\leq n_2\}}\mathbf{1}_{\{z_1> 1\}}m_1(\d \z)\d s\cr\cr
\ar+\ar\frac{1}{k}\int_0^{t\wedge\tau_m}\int_{\mathbb{R}_+^2\setminus\{\mathbf{0}\}}z_1^2\mathbf{1}_{\{\|\z\|
\leq n_2\}}\mathbf{1}_{\{z_1\le 1\}}m_1(\d \z)\d s\cr\cr
\ar+\ar\int_0^{t\wedge\tau_m}\int_{\D_1^c}\zeta_1(s-)^+(e^{z}-1)\mathbf{1}_{\{z\leq n_2\}}\nu(\d z)\d s\cr\cr
 \ar+\ar\int_0^{t\wedge\tau_m}\int_{\D_1}\zeta_1(s-)^+(e^z-1)^2\nu(\d z)\d s + mart.
\eeqnn
Taking expectations on both sides, and letting\ $k\rightarrow\infty,$ we can
find constant\ $C_1$\ large enough such that,
$$\mathbf{E}\zeta_1(t\wedge\tau_m)^+\le
C_1\int_0^t\mathbf{E}[\zeta_1(s\wedge\tau_m)^++\zeta_2(s\wedge\tau_m)^+]\d s.$$
Symmetrically, there exists\ $C_2$ such that,
$$\mathbf{E}\zeta_2(t\wedge\tau_m)^+\le
C_2\int_0^t\mathbf{E}[\zeta_1(s\wedge\tau_m)^++\zeta_2(s\wedge\tau_m)^+]\d s.$$
Let\ $C=C_1+C_2,$
$$\E[\zeta_1(t\wedge\tau_m)^++\zeta_2(t\wedge\tau_m)^+]\le
C\int_0^t\mathbf{E}[\zeta_1(s\wedge\tau_m)^++\zeta_2(s\wedge\tau_m)^+]\d s.$$
By Gronwall's inequality,
$$\E[\zeta_1(t\wedge\tau_m)^++\zeta_2(t\wedge\tau_m)^+]=0, \forall t\ge
0.$$
Since\ $\{\X^{(n_1)}(t):t\ge 0\}$\ and\ $\{\X^{(n_2)}(t):t\ge 0\}$\ are
c\`{a}dl\`{a}g, $\tau_m\rightarrow \infty,$ as $m\rightarrow\infty.$ Thus
$\mathbf{P}\{\X^{(n_1)}(t)\le\X^{(n_2)}(t),\forall t\ge 0\}=1.$
 \qed
 \bigskip

 \bgproposition
$\E(\|\X(t) - \X^{(k)}(t)\|) \rightarrow 0$ and $\X^{(k)}(t) \uparrow \X(t)\ \p-a.s.$
as $k \rightarrow \infty$ for all $t \geq 0.$ Moreover, If\ $\E\|\X(t)\| <
\infty,$ then $\int_1^\infty e^z\nu(\d z) < \infty.$
\edproposition

\proof
Define $\bar{\xi}^{(k)}_i(t):=\xi^{(k)}(t)-b_{ii}t,i=1,2.$ Applying It\^{o}'s formula
to ${\bm e}_1^{\intercal}\X^{(k)}(t)\e^{-\bar{\xi}_1^{(k)}(t)}$ and ${\bm e}_1^{\intercal}\X(t)\e^{-\bar{\xi}^{(k)}_1(t)},$ independently, where ${\bm e}_1=(1,0)$, one can see that
\beqnn
X^{(k)}_1(t)\ar=\ar X_1(0)\e^{\bar{\xi}^{(k)}_1(t)}+\int_0^t\e^{\bar{\xi}^{(k)}_1(t)-\bar{\xi}_1^{(k)}(s)}\sqrt{2c_1X_1^{(k)}(s)}\d B_1(s)\cr\cr
\ar\ar\quad+\int_0^t\int_0^{X^{(k)}_1(s-)}\int_{\mathbb{R}_+^2\setminus\{\mathbf{0}\}}z_1\mathbf{1}_{\{\|\z\|\le k\}}\e^{\bar{\xi}_1^{(k)}(t)-\bar{\xi}_1^{(k)}(s-)}\tilde{M}_1(\d s,\d u,\d \z)\cr\cr
\ar\ar\quad+\int_0^t\int_0^{X^{(k)}_2(s-)}\int_{\mathbb{R}_+^2\setminus\{\mathbf{0}\}}z_1\mathbf{1}_{\{\|\z\|\le k\}}\e^{\bar{\xi}_1^{(k)}(t)-\bar{\xi}_1^{(k)}(s-)}M_2(\d s,\d u,\d \z)\cr\cr
\ar\ar\quad-\int_0^tb_{21}X_2^{(k)}(s)\e^{\bar{\xi}^{(k)}_1(t)-\bar{\xi}^{(k)}_1(s)}\d s
\eeqnn
and
\beqnn
X_1(t)\ar=\ar X_1(0)\e^{\bar{\xi}^{(k)}_1(t)}+\int_0^t\e^{\bar{\xi}^{(k)}_1(t)-\bar{\xi}_1(s)}\sqrt{2c_1X_1(s)}\d B_1(s)\cr\cr
\ar\ar\quad+\int_0^t\int_0^{X_1(s-)}\int_{\mathbb{R}_+^2\setminus\{\mathbf{0}\}}z_1\e^{\bar{\xi}_1^{(k)}(t)-\bar{\xi}_1^{(k)}(s-)}\tilde{M}_1(\d s,\d u,\d \z)\cr\cr
\ar\ar\quad+\int_0^t\int_0^{X_2(s-)}\int_{\mathbb{R}_+^2\setminus\{\mathbf{0}\}}z_1\e^{\bar{\xi}_1^{(k)}(t)-\bar{\xi}_1^{(k)}(s-)}M_2(\d s,\d u,\d \z)\cr\cr
\ar\ar\quad+ \int_0^t \int_{\mathbb{R}} X_1(s-)e^{\bar{\xi}_1^{(k)}(t)-\bar{\xi}_1^{(k)}(s-)} \left(e^{z\mathbf{1}_{\left\{\vert z\vert > k\right\}}} - 1\right) N(\d s, \d z)\cr\cr
\ar\ar\quad-\int_0^tb_{21}X_2(s)\e^{\bar{\xi}^{(k)}_1(t)-\bar{\xi}^{(k)}_1(s)}\d s.
\eeqnn
Hence,
\beqnn
X_1(t)\ar\le\ar X^{(k)}_1(t)+\int_0^t\e^{\bar{\xi}^{(k)}_1(t)-\bar{\xi}^{(k)}_1(s)}\sqrt{2c_1[X_1(s)-X_1^{(k)}(s)]}\d B_1(s)\cr\cr
\ar\ar\quad-\int_0^tb_{21}[X_2(s)-X_2^{(k)}(s)]\e^{\bar{\xi}^{(k)}_1(t)-\bar{\xi}^{(k)}_1(s)}\d s\cr\cr
\ar\ar\quad+\int_0^t\int_{X_1^{(k)}(s-)}^{X_1(s-)}\int_{\mathbb{R}_+^2\setminus\{\mathbf{0}\}}z_1\e^{\bar{\xi}^{(k)}_1(t)-\bar{\xi}^{(k)}_1(s-)}\tilde{M}_1(\d s,\d u,\d \z)\cr\cr
\ar\ar\quad+\int_0^t\int_{X_2^{(k)}(s-)}^{X_2(s-)}\int_{\mathbb{R}_+^2\setminus\{\mathbf{0}\}}z_1\e^{\bar{\xi}^{(k)}_1(t)-\bar{\xi}^{(k)}_1(s-)}M_2(\d s,\d u,\d \z)\cr\cr
\ar\ar\quad+ \int_0^t \int_{\mathbb{R}} X_1(s-)e^{\bar{\xi}^{(k)}_1(t)-\bar{\xi}^{(k)}_1(s-)} \left(e^{z\mathbf{1}_{\left\{\vert z\vert> k\right\}}} - 1\right) N(\d s, \d z).
\eeqnn
Notice that $t \mapsto \E e^{m\bar{\xi}_1^{(k)}(t) }$ is locally bounded for every $m \geq 0$. If $\E \|\X(t)\| < \infty,$ then
$$\left\{\int_0^t e^{-\bar{\xi}^{(k)}_1(s-)}\sqrt{2c_1[X_1(s)-X_1^{(k)}(s)]}\d B_1(s):t\geq0\right\}$$ and
$$\left\{\int_0^t\int_{X_1^{(k)}(s-)}^{X_1(s-)}\int_{\mathbb{R}_+^2\setminus\{\mathbf{0}\}}z_1\e^{-\bar{\xi}^{(k)}_1(s-)}\tilde{M}_1(\d s,\d u,\d \z):t\geq0\right\}$$
 are martingales w.r.t. the filtration $\mcr{F}_t$. Hence
\beqnn
\ar\ar\quad\E\int_0^t e^{\bar{\xi}^{(k)}_1(t)-\bar{\xi}^{(k)}_1(s)}\sqrt{2c_1[X_1(s)-X_1^{(k)}(s)]}\d B_1(s) \cr\cr
\ar=\ar \E\left[e^{\bar{\xi}^{(k)}_1(t)}\E^{\xi}\left(\int_0^t e^{-\bar{\xi}^{(k)}_1(s-)}\sqrt{2c_1[X_1(s)-X_1^{(k)}(s)]}\d B_1(s)\right)\right]\cr\cr
 \ar=\ar 0,
\eeqnn
where $\E^{\xi}$ is the quenched law given $\left\{\xi(t): t\geq 0\right\}$ or $\left\{L(t): t \geq 0\right\}.$ Similarly,
\beqnn
\E\int_0^t\int_{X_1^{(k)}(s-)}^{X_1(s-)}\int_{\mathbb{R}_+^2\setminus\{\mathbf{0}\}}z_1\e^{-\bar{\xi}^{(k)}_1(s-)}\tilde{M}_1(\d s,\d u,\d \z) = 0.
\eeqnn
Taking expectations on both sides,
\beqnn
\ar\ar\E[X_1(t)-X^{(k)}_1(t)]\cr\cr
\ar\le\ar \int_0^t[\int_{\mathbb{R}_+^2}z_1m_2(\d \z)-b_{21}]\E[X_2(s)-X_2^{(k)}(s)]\E\e^{\bar{\xi}^{(k)}_1(t)-\bar{\xi}^{(k)}_1(s)}\d s\cr\cr
\ar\ar\quad+  \int_{\mathbb{R}} \left(e^{z\mathbf{1}_{\left\{\vert z\vert > k\right\}}} - 1\right)\nu (\d z)\int_0^t\E X_1(s-)\E e^{\bar{\xi}^{(k)}_1(t)-\bar{\xi}^{(k)}_1(s-)}  \d s.
\eeqnn
Symmetrically,
\beqnn
\ar\ar\E[X_2(t)-X^{(k)}_2(t)]\cr\cr
\ar\le\ar \int_0^t[\int_{\mathbb{R}_+^2}z_2m_1(\d \z)-b_{12}]\E[X_1(s)-X_1^{(k)}(s)]\E\e^{\bar{\xi}^{(k)}_1(t)-\bar{\xi}^{(k)}_1(s)}\d s\cr\cr
\ar\ar\quad+  \int_{\mathbb{R}} \left(e^{z\mathbf{1}_{\left\{\vert z\vert > k\right\}}} - 1\right)\nu (\d z)\int_0^t\E X_2(s-)\E e^{\bar{\xi}^{(k)}_1(t)-\bar{\xi}^{(k)}_1(s-)}  \d s.
\eeqnn
Since  $\E\|\X(t)\| <\infty,$ $\int_{\mathbb{R}} \left(e^{z\mathbf{1}_{\left\{\vert z\vert > k\right\}}} - 1\right)\nu (\d z)<\infty.$ Thus,
$$\int_{\mathbb{R}} \left(e^{z\mathbf{1}_{\left\{\vert z\vert > k\right\}}} - 1\right)\nu (\d z)=\int_{\mathbb{R}} \left(e^{z} - 1\right)\mathbf{1}_{\left\{\vert z\vert > k\right\}}\nu (\d z)$$
tends to $0$ as $k\rightarrow\infty$. Then by taking limits of $k\rightarrow\infty$ and Gronwall's inequality, we get the conclusion.
\qed

\bigskip\bigskip

\noindent{\it Proof of Theorem 1.1:}

\textbf{Step 1.}
Denote $\tilde{\beta} = a + \frac{1}{2}\sigma_1^2 + \int_{\D_1}(e^{z} - 1 - z)\nu(\d z) + \int_{\D_1^c}(e^{z} - 1)\nu(\d z), $
By Lemma 3.2 in \cite{ref7},
$$
\E \e^{\xi(t)}=\e^{\tilde{\beta}t},\ \forall t\ge 0.
$$
According to the formula (40) in \cite{qinzheng},
$$
\E^{\xi}\Big[{\mrm{exp}}\{-\langle{\bm \lambda},\X(t)\rangle\}\vert\mathcal{F}_r\Big]=\exp\{-\langle \X(r),{\bm v}_{r,t}(\xi,{\bm\lambda})\rangle\},\,\,t\ge r\ge0.
$$
Taking derivatives with respect to\ ${\bm\lambda}=\mathbf{0}+$\ on both sides, we have, $$\E^{\xi}(\X(t)\vert\mathcal{F}_r)=\\e^{\xi(t)-\xi(r)}\e^{(r-t)\tilde{\b}^{\intercal}}\X(r).$$
Hence for any\ $\mathcal{F}_r$-measurable function\ $F,$
\beqnn
  \E\big(F\e^{-\tilde{\beta}t}\e^{t\tilde{\b}^{\intercal}}\X(t)\big)\ar=\ar
  \E\big(F\e^{-\tilde{\beta}t}\e^{t\tilde{\b}^{\intercal}}
  \E^{\xi}(\X(t)\vert\mathcal{F}_r)\big)\cr
\ar=\ar\e^{-\tilde{\beta}t}\e^{t\tilde{\b}^{\intercal}}\E\big(F\e^{\xi(t)-\xi(r)}\e^{(r-t)\tilde{
\b}^{\intercal}}\X(r)\big)\cr
\ar=\ar\e^{-\tilde{\beta}t}\e^{t\tilde{\b}^{\intercal}}\E[\e^{\xi(t)-\xi(r)}]\E\big(F\e^{(r-t)\tilde{
\b}^{\intercal}}\X(r)\big)\cr
\ar=\ar\e^{-\tilde{\beta}t}\e^{t\tilde{\b}^{\intercal}}\e^{\tilde{\beta}t-\tilde{\beta}r}\E
\big(F\e^{(r-t)\tilde{\b}^{\intercal}}\X(r)\big)\cr
\ar=\ar\E\big(F\e^{-\tilde{\beta}r}\e^{r\tilde{\b}^{\intercal}}\X(r)\big),
\eeqnn
where $\tilde{\b}$\ is a $2\times2$-matrix with $\tilde{b}_{11}=b_{11},$~$\tilde{b}_{22}=b_{22},$~$\tilde{b}_{12}=b_{12}-\int_{\mathbb{R}_+^2}z_2m_1(\d \z),$~$\tilde{b}_{21}=b_{21}-\int_{\mathbb{R}_+^2}z_1m_2(\d\z).$
Thus,\ ${\bm{\mathcal{M}}}(t)=\e^{-\tilde{\beta}t}\e^{t\tilde{\b}^{\intercal}}\X(t)$\ is a two-dimensional martingale, and
\beqnn
  \E\X(t)=\Big[\E\frac{{\bm v}_{0,t}({\bm\xi},{\bm\lambda})}{\d {\bm\lambda}}\big\vert_{{\bm\lambda}={\bm0}+}\Big]^{\intercal}\E\X(0) =\e^{\tilde{\beta}t}\E \e^{-t\tilde{\b}^{\intercal}}\X(0).
\eeqnn

\textbf{Step 2.}
Denote $$\beta(n) = an + \frac{\sigma^2n^2}{2} + \int_{\D_1}(e^{nz} - 1 - nz)\nu(\d z) + \int_{\D_1^c}(e^{nz} - 1)\nu(\d z).$$
Since $\int_1^{\infty}\e^{nz}\nu(\d z)<\infty,$\ by similar arguments of Lemma 3.2 in \cite{ref7}, for any integer\ $m\le n,$ we have
\begin{equation}\label{2020020301}
\E\e^{m\bar{\xi}_1(s)}\le C_m(t):=\exp\{[(\beta(m)-b_{11}m)t]\vee 0\}<\infty,\ s\in[0,t],
\end{equation}
where $\bar{\xi}_1(s)=\xi(s)-b_{11}s.$ Applying It\^{o}'s formula to\ $[X_1(t)\e^{-\bar{\xi}_1(t)}]^n,$

\begin{eqnarray*}
\ar\ar[X_1(t)\e^{-\bar{\xi}_1(t)}]^n\cr\cr
\ar=\ar \sum\limits_{j=0}^{n-2}A^1_{n,j}\int_0^t[X_1(s)]^{j+1}\e^{-n\bar{\xi}_1(s)}\d s -\int_0^tb_{21}n[X_1(s)]^{n-1}X_2(s)\e^{-n\bar{\xi}_1(s)}\d s\cr \ar\quad\ar+ n\int_0^t [X_1(s)]^{n-1}\e^{-n\bar{\xi}_1(s)}\sqrt{2c_1X_1(s)}\d B_1(s)+[X_1(0)]^n\cr\cr
\ar\quad\ar+ \sum_{j = 0}^{n - 1}\binom{n}{j} \int_0^t\int_0^{X_1(s-)}\int_{\mathbb{R}_+^2\setminus\{\mathbf{0}\}}X_1(s-)^j z_1^{n - j}e^{-n\bar{\xi}_1(s-)}\tilde{M_1}(\d s, \d u, \d\z)\cr\cr
\ar\quad\ar+\int_0^t\int_0^{X_2(s-)}\int_{\mathbb{R}_+^2\setminus\{\mathbf{0}\}}\big([X_1(s)+z_1]^n-[X_1(s)]^n\big)\e^{-n\bar{\xi}_1(s)}M_2(\d s,\d u,\d\z),
\end{eqnarray*}
By a standard stopping time argument, we can see that
$$\int_0^t [X_1(s)]^{n-1}\e^{-n\bar{\xi}_1(s)}\sqrt{2c_1X_1(s)}\d B_1(s)$$
and
$$\int_0^t\int_0^{X_1(s-)}\int_{\mathbb{R}_+^2\setminus\{\mathbf{0}\}}X_1(s-)^j z_1^{n - j}e^{-n\bar{\xi}_1(s-)}\tilde{M_1}(\d s, \d u, \d\z)$$
are local martingales. Therefore, there exists a sequence of stopping times $\tau_k$ such that
\beqnn
\ar\ar\E^{\xi} [X_1(t\wedge\tau_k)\e^{-\bar{\xi}_1(t\wedge\tau_k)}]^n
=\E^{\xi}[X_1(0)]^n +\sum\limits_{j=0}^{n-2}A^1_{n,j}\E^{\xi}\int_0^{t\wedge\tau_k}[X_1(s)]^{j+1}\e^{-n\bar{\xi}_1(s)}\d s\cr
\ar\ar\qquad\qquad\qquad\qquad\qquad +\sum\limits_{j=0}^{n-1}B^1_{n,j}\E^{\xi}\int_0^{t\wedge\tau_k}[X_1(s)]^{j}X_2(s)\e^{-n\bar{\xi}_1(s)}\d s,
\eeqnn
Therefore,
\begin{eqnarray}\label{2020020304}
\ar\ar\E[X_1(t\wedge\tau_k)]^n\cr\cr
\ar\ar=\E[X_1(0)]^n \E\e^{n\bar{\xi}_1(t\wedge\tau_k)} +\sum\limits_{j=0}^{n-2}A^1_{n,j}\E\big(\e^{n\bar{\xi}_1(t\wedge\tau_k)}\int_0^{t\wedge\tau_k}[X_1(s)]^{j+1}\e^{-n\bar{\xi}_1(s)}\d s\big)\cr
\ar\ar\qquad\qquad\qquad +\sum\limits_{j=0}^{n-1}B^1_{n,j}\E\big(\e^{n\bar{\xi}_1(t\wedge\tau_k)}\int_0^{t\wedge\tau_k}[X_1(s)]^{j}X_2(s)\e^{-n\bar{\xi}_1(s)}\d s\big).
\end{eqnarray}
By Fubini's Theorem,
\begin{eqnarray}\label{2020020302}
\ar\ar\E\big(\e^{n\bar{\xi}_1(t\wedge\tau_k)}\int_0^{t\wedge\tau_k}[X_1(s)]^{j}X_2(s)\e^{-n\bar{\xi}_1(s)}\d s\big)\cr
\ar\ar\qquad\qquad\qquad\le  \int_0^t\E [X_1(s\wedge\tau_k)]^{j}\E X_2(s\wedge\tau_k)\E\e^{n[\bar{\xi}_1(t\wedge\tau_k)-\bar{\xi}_1(s\wedge\tau_k)]}
\end{eqnarray}
and
\begin{eqnarray}\label{2020020303}
\ar\ar\E\big(\e^{n\bar{\xi}_1(t\wedge\tau_k)}\int_0^{t\wedge\tau_k}[X_1(s)]^{j+1}\E^{-n\bar{\xi}_1(s)}\d s\big)\cr
\ar\ar\qquad\qquad\qquad\le  \int_0^t\E [X_1(s\wedge\tau_k)]^{j+1}\E\e^{n[\bar{\xi}_1(t\wedge\tau_k)-\bar{\xi}_1(s\wedge\tau_k)]}.
\end{eqnarray}
Substituting \eqref{2020020301},\ \eqref{2020020302},\ \eqref{2020020303} into \eqref{2020020304} we get,
\beqnn
\E [X_1(t\wedge\tau_k)]^n\ar\le\ar C_n(t)\bigg[\E [X_1(0)]^n + \sum\limits_{j=0}^{n-2}A^1_{n,j}\int_0^t\E [X_1(s\wedge\tau_k)]^{j+1}\d s\cr
\ar\ar\qquad+\sum\limits_{j=0}^{n-1}B^1_{n,j}\int_0^t\E[X_1(s\wedge\tau_k)]^{j}\E X_2(s\wedge\tau_k)\d s\bigg].
\eeqnn
By step 1, ${\bm{\mathcal{M}}}(t)=\e^{-\tilde{\beta}t}\e^{t\tilde{\b}^{\intercal}}\X(t)$\ is a martingale. Hence, for any positive integer $m$ and $s\in[0,t],$
\begin{eqnarray}\label{2020020401}
\E [X_1(s)]^m\ar=\ar\e^{m\tilde{\beta}s}\E [{e}^{\intercal}_1\e^{-s\tilde{\b}^{\intercal}}\bm{\mathcal{M}}(s)]^m\cr
\ar\le\ar\e^{m\tilde{\beta}s}\E [\|\e^{-s\tilde{\b}^{\intercal}}\bm{\mathcal{M}}(s)\|^m]\cr
\ar\le\ar\e^{m\tilde{\beta}s}\|\e^{-s\tilde{\b}^{\intercal}}\|^m\E [\|\bm{\mathcal{M}}(s)\|^m]\cr
\ar\le\ar\e^{m\tilde{\beta}s}\|\e^{-s\tilde{\b}^{\intercal}}\|^m(\sqrt{2})^m\E [(\mathcal{M}_1(s))^m+(\mathcal{M}_2(s))^m]\cr
\ar\le\ar\e^{m\tilde{\beta}s}\|\e^{-s\tilde{\b}^{\intercal}}\|^m(\sqrt{2})^m\E [(\mathcal{M}_1(t))^m+(\mathcal{M}_2(t))^m]\cr
\ar\le\ar\e^{m\tilde{\beta}s}\|^m\e^{-s\tilde{\b}^{\intercal}}\|^m2(\sqrt{2})^m\E \|\bm{\mathcal{M}}(t)\|^m\cr
\ar\le\ar\e^{m\tilde{\beta}(s-t)}\|\e^{-s\tilde{\b}^{\intercal}}\|^m\|\e^{t\tilde{\b}^{\intercal}}\|^m2(\sqrt{2})^m\E \|\X(t)\|^m\cr
\ar\le\ar O(t,m)\E \|\X(t)\|^m,
\end{eqnarray}
where $O(t,m):=\sup\limits_{n\le m}\sup\limits_{s\in[0,t]}\e^{n\tilde{\beta}(s-t)}\|\e^{-s\tilde{\b}^{\intercal}}\|^n\|\e^{t\tilde{ \b}^{\intercal}}\|^n2(\sqrt{2})^n.$ Symmetrically,
\beqnn
\E [X_2(t\wedge\tau_k)]^n\ar\le\ar C'_n(t)\bigg[\E [X_2(0)]^n + \sum\limits_{j=0}^{n-2}A^2_{n,j}\int_0^t\E [X_2(s\wedge\tau_k)]^{j+1}\d s\cr
\ar\ar\qquad+\sum\limits_{j=0}^{n-1}B^2_{n,j}\int_0^t\E[X_2(s\wedge\tau_k)]^{j}\E X_1(s\wedge\tau_k)\d s\bigg]
\eeqnn
and
\begin{equation}\label{2020020403}
\E [X_2(s)]^m\le O(t,m)\E \|\X(t)\|^m, \quad s\in[0,t],
\end{equation}
where $C'_m(t):=\exp\{[(\beta(m)-b_{22}m)t]\vee 0\}. $ Define\ $\bar{C}_n(t):=C_n(t)\vee C'_n(t),$ $A_{n,j}:=A_{n,j}^1\vee A_{n,j}^2,$ $B_{n,j}:=B_{n,j}^1\vee B_{n,j}^2.$
From above, it is obvious that for $i=1,2,$
\beqnn
\E [X_i(t\wedge\tau_k)]^n\ar\le\ar \bar{C}_n(t)\bigg[\E [X_i(0)]^n + t\sum\limits_{j=0}^{n-2}A_{n,j}O(t,n)\E \|\X(t)\|^{j+1}\d s\cr
\ar\ar\qquad+t\sum\limits_{j=0}^{n-1}B_{n,j}O(t,n)\E\|\X(t)\|^{j+1}\d s\bigg].
\eeqnn
By \eqref{2020020401} and \eqref{2020020403}, $t\mapsto\E\|\X(t)\|^m$ is locally bounded under the asumption that $\E\|\X(t)\|^m<\infty.$ By induction one can find a finite function $g_n(t)$
s.t.
\beqnn
g_n(t)\ar\ge\ar \bar{C}_n(t)\bigg[\E [X_i(0)]^n + t\sum\limits_{j=0}^{n-2}A_{n,j}O(t,n)\E \|\X(t)\|^{j+1}\d s\cr
\ar\ar\qquad+t\sum\limits_{j=0}^{n-1}B_{n,j}O(t,n)\E\|\X(t)\|^{j+1}\d s\bigg].
\eeqnn
Thus, by Fatou's Lemma, for $i=1,2,$
$$\E[X_i(t)]^n\le\liminf\limits_{k\rightarrow\infty}\E [X_i(t\wedge\tau_k)]^n\le g_n(t)<\infty,$$
which implies $\E\|\X(t)\|^n<\infty.$

\textbf{Step 3.}
Applying It\^{o}'s formula to $[{\bm e}_1^{\intercal}\X^{(k)}(t)\exp\{-\bar{\xi}^{(k)}_1(t)\}]^n,$
\begin{eqnarray*}
\ar\ar\quad [X_1^{(k)}(t)\e^{-\bar{\xi}^{(k)}_1(t)}]^n\cr\cr
\ar=\ar [X_1(0)]^n -\int_0^tb_{21}n[X^{(k)}_1(s)]^{n-1}X^{(k)}_2(s)\e^{-n\bar{\xi}^{(k)}_1(s)}\d s\cr\cr
 \ar\quad\ar+ n\int_0^t [X^{(k)}_1(s)]^{n-1}\e^{-n\bar{\xi}^{(k)}_1(s)}\sqrt{2c_1X^{(k)}_1(s)}\d B_1(s)\cr\cr
\ar\quad\ar +\int_0^t n(n-1)c_1[X^{(k)}_1(s)]^{n-1}e^{-n\bar{\xi}^{(k)}_1(s-)}\d s\cr\cr
\ar\quad\ar-\int_0^t\int_{\mathbb{R}_+^2\setminus\{\mathbf{0}\}}n [X^{(k)}_1(s)]^{n}\e^{-n\bar{\xi}^{(k)}_1(s)}z_1\mathbf{1}_{\{\|\z\|\le k\}}m_1(\d \z)\d s\cr\cr
\ar\quad\ar+ \sum_{j = 0}^{n - 1}\binom{n}{j} \int_0^t\int_0^{X^{(k)}_1(s-)}\int_{\mathbb{R}_+^2\setminus\{\mathbf{0}\}}X^{(k)}_1(s-)^j z_1^{n - j}\mathbf{1}_{\{\|\z\|\le k\}}\e^{-n\bar{\xi}^{(k)}_1(s-)}{M_1}(\d s,\d u,\d\z)\cr\cr
\ar\quad\ar+ \sum_{j = 0}^{n - 1}\binom{n}{j} \int_0^t\int_0^{X^{(k)}_2(s-)}\int_{\mathbb{R}_+^2\setminus\{\mathbf{0}\}}X^{(k)}_1(s-)^j z_1^{n - j}\mathbf{1}_{\{\|\z\|\le k\}}\e^{-n\bar{\xi}^{(k)}_1(s-)}M_2(\d s,\d u,\d\z).
\end{eqnarray*}
 Since $\E \|\X^{(k)}(t)\|<\infty$ and $\E \e^{n\bar{\xi}^{(k)}_1(t)}<\infty$ for each positive integer $n,$ and $\X^{(k)}(t) \uparrow\X(t),$ $\bar{\xi}^{(k)}_1(t) \uparrow \bar{\xi}_1(t),$ almost surely for $\p$ as $k\rightarrow\infty$, we have
\begin{eqnarray*}
\ar\ar[X^{(k)}_1(t)]^n\cr\cr
\ar=\ar [X^{(k)}_1(0)]^n\e^{n\bar{\xi}_1(t)} -\int_0^tb_{21}n[X^{(k)}_1(s)]^{n-1}X^{(k)}_2(s)\e^{n(\bar{\xi}^{(k)}_1(t)-\bar{\xi}^{(k)}_1(s))}\d s\cr\cr
\ar\quad\ar +\int_0^t n(n-1)c_1[X^{(k)}_1(s)]^{n-1}e^{n\bar{\xi}^{(k)}_1(t)-n\bar{\xi}_1^{(k)}(s-)}\d s\cr\cr
\ar\quad\ar+ \sum_{j = 0}^{n - 2}\binom{n}{j} \int_0^t\int_{\mathbb{R}_+^2\setminus\{\mathbf{0}\}}X^{(k)}_1(s-)^{j+1} z_1^{n - j}\mathbf{1}_{\{\|\z\|\le k\}}\e^{n\bar{\xi}^{(k)}_1(t)-n\bar{\xi}^{(k)}_1(s-)}m_1(\d\z)\d s\cr\cr
\ar\quad\ar+ \sum_{j = 0}^{n - 1}\binom{n}{j} \int_0^t\int_{\mathbb{R}_+^2\setminus\{\mathbf{0}\}}X^{(k)}_1(s-)^j X^{(k)}_2(s-) z_1^{n - j}\mathbf{1}_{\{\|\z\|\le k\}}\e^{n\bar{\xi}^{(k)}_1(t)-n\bar{\xi}^{(k)}_1(s-)}m_2(\d \z)\d s\cr\cr
\ar\quad\ar+mart.
\end{eqnarray*}
Taking expectation on both sides and letting $k\rightarrow\infty,$ we get the desired result by monotone convergence theorem.
\qed
\section{$f$-moment of two-type CBRE processes}

\setcounter{equation}{0}

In this section, we present the equivalent condition for the existence of $\E f(\|\X(t)\|),$ where
$f$ is a nonnegative continuous function $[0,\infty)$ satisfying \textbf{Condition A.}
By discussions in \cite{ath}, pp.154, this condition can be changed by

\noindent\textbf{Condition B.} There exists constants $K>0$ such that

\quad\quad\textbf{(B1)} $f$ is convex and nondecreasing on $[0,\infty);$

\quad\quad\textbf{(B2)} For all  $x,y\in[0,\infty),$ $f(xy)\le Kf(x)f(y)$;

\quad\quad\textbf{(B3)} For all $x\in[0,\infty),$ $f(x)>1.$

\bigskip

\bgproposition\label{prop2020012901}
 Suppose that $f$ satisfies \textbf{Condition B}. Then for any $t\ge 0$ and $\y\ge \x\in\mathbb{R}_+^2,$
 $\E f(\|\X(t,\y)\|)<\infty$ if and only if $\E f(\|\X(t,\x)\|)<\infty,$ where $\{\X(t,\x):t\geq0\}$ is a unique strong solution on $\mbb{R}_{+}^2$ to (\ref{eq01})--(\ref{eq02}) starting from any fixed point $\x\in\mbb{R}_+^2$.
\edproposition

\proof
Denote $\X^k(t,\x)=\X(t,k\x)-\X(t,(k-1)\x).$ By the quenched branching property, $\{\X^k(t,\x):t\ge 0\},\ k=1,2,....$ are independently identically distributed. For $z\in[0,\infty),$ denote $\lfloor z\rfloor$ as the integer part of $z.$ For any $\y\ge \x\in\mathbb{R}_+^2,$ define $r({\x,\y}):=\lfloor y_1/x_1\rfloor\vee\lfloor y_2/x_2\rfloor.$ It is clear that $r(\x,\y)\le\|(y_1/x_1,y_2/x_2)\|.$ By \textbf{Condition B},
\beqnn
\ar\ar\E^{\xi}[f(\|\X(t,\y)\|)\mathbf{1}_{\{t<\tau_n(\y)\}}]\cr
\ar=\ar\E^{\xi}\big[f\big(\|\sum\limits_{k=1}^{r(\x,\y)}\X^k(t,\x)+\X(t,\y)-\X(t,r(\x,\y)\x)\|\big)\mathbf{1}_{\{t<\tau_n(\y)\}}\big]\cr\cr
\ar\le\ar Kf(1+\|(y_1/x_1,y_2/x_2)\|)\E^{\xi}[f(\|\X(t,\x)\|)\mathbf{1}_{\{t<\tau_n(\x)\}}],
\eeqnn
where $\tau_0(\x)=0,$\ $\tau_n(\x)$ represents the $n^{\text{th}}$ jumping time that the jump size of $\{\X(t,\x):t\ge0\}$ falls into $\D_2:=\mathbb{R}^2_+\setminus[0,1]^2.$ Taking expectation on both sides and letting $n\rightarrow\infty,$
\begin{equation}\label{0708}
  \E[f(\|\X(t,\y)\|)]\le Kf(1+\|(y_1/x_1,y_2/x_2)\|)\E [f(\|\X(t,\x)\|)].
\end{equation}
The desired result follows.\qed

\bgproposition\label{prop2020012902}
 Suppose that $f$ satisfies \textbf{Condition B}, and $\E f(\|\X(t,\x)\|)<\infty$ for some $\x\in\mathbb{R}_+^2$ and some $ t\ge 0.$ Then $\E f(\|\X(t)\|)<\infty$ if and only if $\E f(\|\X(0)\|)<\infty.$
\edproposition
\proof
Some simple calculations lead to,
$$
  \E f(\|\X(t)\|)\le \frac{1}{2}K^2f(2)[f(1)+\E f(\|\X(0)\|)]\E [f(\|\X(t,\mathbf{1})\|)].
$$
By Proposition \ref{prop2020012901} we have
 $\E f(\|\X(t,\mathbf{1})\|)<\infty.$ Then $\E f(\|\X(0)\|)<\infty$ leads to $\E f(\|\X(t)\|)<\infty.$ Conversely, we suppose that $\E f(\|\X(t)\|)<\infty.$ According to Step 1, $\bm{\mathcal{M}}(t)=\e^{-\tilde{\beta}t}\e^{\tilde{\bf b}t}\X(t)$ is a two-dimensional martingale, then
$$\E f(\|\bm{\mathcal{M}}(t)\|)\le K^2f(\e^{\tilde{\beta}t})f(\|\e^{t\tilde{\b}}\|)\E f(\|\X(t)\|)<\infty.$$
Moreover, by the convexity of $f$,
\beqnn
\E f(\|\X(0)\|)\ar=\ar\E f(\|\bm{\mathcal{M}}(0)\|)\le\E f(\mathcal{M}_1(0)+\mathcal{M}_2(0))\cr
\ar\le\ar\frac{1}{2}Kf(2)[\E f(\mathcal{M}_1(0))+\E f(\mathcal{M}_2(0))]\cr
\ar\le\ar\frac{1}{2}Kf(2)[\E f(\E[\mathcal{M}_1(t)\vert\mathcal{F}_0])+\E f(\E[\mathcal{M}_2(t)\vert\mathcal{F}_0])]\cr
\ar\le\ar Kf(2)\E f(\E[\|\bm{\mathcal{M}}(t)\|\vert\mathcal{F}_0])\cr
\ar\le\ar Kf(2)\E (\E[f(\|\bm{\mathcal{M}}(t)\|)\vert\mathcal{F}_0])\cr
\ar=\ar Kf(2)\E f(\|\bm{\mathcal{M}}(t)\|).
\eeqnn
Finally, the desired result follows.\qed

\bglemma\label{lem2020020401}
 Suppose that $f$ satisfies \textbf{Condition B}, and for any $n\ge 1,$ $\int_{\{\|\z\|\ge 1\}}\|\z\|^n(m_1+m_2)(\d\z)$ and $\int_1^{\infty}\e^{nz}\nu(\d z)$ are finite. Then for any $\x\in\mathbb{R}_+^2,$ $t\mapsto\E f(\|\X(t,\x)\|)$ is locally bounded on $[0,\infty).$
\edlemma

\proof
 From the previous section, we have $\E  f(\|\X(t,\x)\|^n )<\infty.$ In view of pp.160 in \cite{Sato} and the proof of Lemma 3.7 in \cite{ref6}, there exist a constant $c>0$ and a positive integer $n$ such that for all $z\ge 1,$ $f(\e^{nz})\le c \e^{nz}.$ Then
$$\E f(\|\X(t,\x)\|)\le f(1)+c\E[\|\X(t,\x)\|^n]<\infty.$$
Since for $i=1,2,$ $f(\mathcal{M}_i(t,\x))$ is an $\mathcal{F}_t$-submartingale, for $t\in[0,T],$
\beqnn
\E f(\|\bm{\mathcal{M}}(t,\x)\|)\ar\le\ar \E f(\mathcal{M}_1(t,\x)+\mathcal{M}_2(t,\x))\cr
\ar\le\ar \frac{1}{2}Kf(2)\big[\E f(\mathcal{M}_1(t,\x))+\E f(\mathcal{M}_2(t,\x))\big]\cr
\ar\le\ar\frac{1}{2}Kf(2)\big[\E f(\mathcal{M}_1(T,\x))+\E f(\mathcal{M}_2(T,\x))\big]\cr
\ar\le\ar Kf(2)\E f(\|\bm{\mathcal{M}}(T,\x)\|)\cr
\ar\le\ar K^3f(2)\E f(\|\e^{-\tilde{\beta}T}\|)f(\|\e^{T\tilde{\b}}\|\vee 1)f(\|\X(T,\x)\|)<\infty.
\eeqnn
Hence,
\begin{eqnarray*}
\E f[\|\X(t,\x)\|]\ar=\ar\E f[\|\e^{\tilde{\beta}t}\e^{-t\tilde{\b}}\e^{\tilde{\beta}t}\e^{t\tilde{\b}}\X(t,\x)\|]\le K f(\|\e^{\tilde{\beta}t}\e^{-t\tilde{\b}}\|)\E f(\|\bm{\mathcal{M}}(t,\x)\|)\cr
\ar\le\ar K^4 f(2) f(\|\e^{\tilde{\beta}t}\e^{-t\tilde{\b}}\|) f(\|\e^{-\tilde{\beta}T}\|)f(\|\e^{T\tilde{\b}}\|\vee 1)\E (\|\X(T,\x)\|).
\end{eqnarray*}
Therefore $t\mapsto\E f(\|\X(t,\x)\|)$ is locally bounded on $[0,\infty).$
\qed

For $\x\in\mathbb{R}^2_+,$ let $\theta_0(\x)=0$ and $\theta_n(\x) = \theta_n'(\x) \wedge \theta_n''(\x)$ for $n\geq1$, where
\beqnn
&&\theta_n'(\x) = \inf \left\{t > \theta_{n -1}(\x): \xi_t - \xi_{t-} > 1\right\},\\
&&\theta_n''(\x) = \inf \left\{t > \theta_{n - 1}(\x): [X_1(t,\x) - X_1(t-,\x)]\wedge [X_2(t,\x) - X_2(t-,\x)]> 1, \xi_t = \xi_{t-}\right\}.
\eeqnn
Let $J(\d t)$ be the distribution of $\theta_1(\mathbf 1),$ that is $J(\d t)=\E(\theta_1(\mathbf 1)\in \d t).$ Define
$$
\mu_n(t):=\E [f(\|\X(t,\mathbf{1})\|);t<\theta_1(\mathbf 1)].
$$
Notice that $\mu_0(t) = 0.$

\bgproposition\label{prop2020020402}
 Suppose that $f$ satisfies \textbf{Condition B} and
 $$
 \int_{\{\|\z\|\ge 1\}}f(\|\z\|)(m_1+m_2)(\d \z)<\infty,\,\,\int_1^{\infty}f(\e^z)\nu(\d z)<\infty.
 $$
 Then, for any $\x\in\mathbb{R}_+^2,$ $t\mapsto\E f(\|\X(t,\x)\|)$ is locally bounded on $[0,\infty).$
\edproposition
\proof
Recall that $\{\X(t,\bm{1}):t\ge  0\}$ is the strong solution with initial value $\mathbf{1}=(1,1).$ On the same probability space, define ${\bm R}(t,\x)$ to the strong solution of the following equation system:
\begin{eqnarray}\label{2020020406}
R_1(t)&=&x_1 + c_1\int_0^t\int_0^{R_1(s)}W_1(\d s,\d u) + \int_0^t(-b_{11}R_1(s)-b_{21}R_2(s))\d s\cr\cr
& &+\int_0^t\int_0^{R_1(s-)}\int_{\D_2^c}z_1\tilde{M}_1(\d s,\d u,\d\z)\cr\cr
& &+\int_0^t\int_0^{R_2(s-)}\int_{\D_2^c}z_1M_2(\d s,\d u,\d\z)\cr\cr
& &+\int_0^t\int_{\D_1}R_1(s-)(e^z-1)\tilde{N}(\d s,\d z)\cr\cr
& &+\int_0^t\int_{-\infty}^{-1}R_1(s-)(e^z-1)N(\d s,\d z)
\end{eqnarray}
and
\begin{eqnarray}\label{2020020407}
R_2(t)&=&x_2 +  c_2\int_0^t\int_0^{R_2(s)}W_2(\d s,\d u)  + \int_0^t(-b_{12}R_1(s)-b_{22}R_2(s))\d s\cr\cr
& &+\int_0^t\int_0^{R_1(s-)}\int_{\D_2^c}z_2M_1(\d s,\d u,\d\z)\cr\cr
& &+\int_0^t\int_0^{R_2(s-)}\int_{\D_2^c}z_2\tilde{M}_2(\d s,\d u,\d\z)\cr\cr
& &+\int_0^t\int_{\D_1}R_2(s-)(e^z-1)\tilde{N}(\d s,\d z)\cr\cr
& &+\int_0^t\int_{-\infty}^{-1}R_2(s-)(e^z-1)N(\d s,\d z).
\end{eqnarray}
 Let $W'$ be the space consisting of all c\`{a}dl\`{a}g paths $t\mapsto\x(t)$ from $[0,\infty)$ to $\mathbb{R}^2_+$ with Skorokhod topoligy. Let $\mathfrak{G}=\sigma\{\x(s):s\ge 0\},$ $\mathfrak{G}_t=\sigma\{\x(s):0\le s\le t\},$\ $t\ge 0$ be natural filtrations on $W'.$ Denote $\mbb{P}_{\x}$ the distribution of $\{\X(t,\x):t\ge 0\}$ on $W',$ then $(W',\mathfrak{G},\mathfrak{G}_t,\mbb{P}_{\x})$ is a canonical realization of the two-dimensional CBRE process with branching mechanism $\bm{\phi}.$ Denote $\sigma_n$ the stopping time of $\{\x(t) : t \ge 0\}$ corresponding to the stopping time
$\theta_n(x)$ of $\{\X(t,\x) : t \ge 0\}$. And $\mbb{E}_x$ stands for the mathematical expectation with repest to $\mbb{P}_x.$ Then,
\begin{eqnarray*}
\mu_n(t)\ar=\ar\E[f(\|\X(t,\mathbf{1})\|)\mathbf{1}_{\{t<\theta_1(\mathbf{1})\}}]+\E[f(\|\X(t,\mathbf{1})\|)\mathbf{1}_{\{\theta_1(\mathbf{1})\le t<\theta_n(\mathbf{1})\}}]\cr
\ar=\ar \E[f(\|{\bm R}(t,\mathbf{1})\|)+\E\{\mathbf{1}_{\{\theta_1(\mathbf{1})\le t\}}\E[f(\|\X(t,\mathbf{1})\|)\mathbf{1}_{\{ t<\theta_n(\mathbf{1})\}}\vert\mathfrak{G}_{\theta_1(\mathbf{1})}]\}\cr
\ar\le\ar \E[f(\|{\bm R}(t,\mathbf{1})\|)+\E\{\mathbf{1}_{\{\theta_1(\mathbf{1})\le t\}}\mbb{E}_{\X(\theta_1(\mathbf{1}))}[f(\|\x(t-\theta_1(\mathbf{1}))\|)\mathbf{1}_{\{t-\theta_1(\mathbf{1})<\sigma_{n-1}\}}]\}]\cr
\ar\le\ar \E[f(\|{\bm R}(t,\mathbf{1})\|)\cr
\ar\ar\qquad+\E\{\mathbf{1}_{\{\theta_1(\mathbf{1})\le t\}}\mbb{E}_{({\bm R}(\theta_1(\mathbf{1}))+ \Delta\X(\theta_1(\mathbf{1})))}[f(\|\x(t-\theta_1(\mathbf{1})\|)\mathbf{1}_{\{t-\theta_1(\mathbf{1})<\sigma_{n-1}\}}]\}].
\end{eqnarray*}
Without loss of generalization, suppose $m_1(\D_2),$\ $m_2(\D_2)$ and $\nu(1,\infty)$ are positive. Denote
$$\hat{m}_1(\d \z)=\frac{\mathbf{1}_{\{\z\in\D_2\}}}{m_1(\D_2)}m_1(\d \z),$$
$$\hat{m}_2(\d\z)=\frac{\mathbf{1}_{\{\z\in\D_2\}}}{m_2(\D_2)}m_2(\d \z),$$
$$\hat{\nu}(\d z)=\frac{\mathbf{1}_{\{z\in(1,\infty)\}}}{\nu(1,\infty)}\nu(\d z).$$
 By assumption, $M_1(\d s,\d u,\d\z),$ $M_2(\d s,\d u,\d\z)$ and $N(\d s,\d z)$ are mutually independent,
\begin{eqnarray*}
\ar\ar\qquad\E\{\mathbf{1}_{\{\theta_1(\mathbf{1})\le t\}}\mbb{E}_{({\bm R}(\theta_1(\mathbf{1}))+\Delta \X(\theta_1(\mathbf{1})))}[f(\|\x(t-\theta_1(\mathbf{1}))\|)\mathbf{1}_{\{t-\theta_1(\mathbf{1})<\sigma_{n-1}\}}]\}\cr\cr
\ar\le\ar\int_0^tJ(\d
s)\int_{\D_2}\int_{\D_2}\int_1^{\infty}\E\{\mbb{E}_{(\e^z
{\bm R}(s,\mathbf{1})+{\bm u}+{\bm v})}[f(\|\x(t-s)\|)
\mathbf{1}_{\{t-s<\sigma_{n-1}\}}]\}
\hat{m}_1(\d {\bm u})\hat{m}_2(\d{\bm v})\hat{\nu}(\d z).
\end{eqnarray*}
According to \eqref{0708}, the above is no more than
\begin{eqnarray*}
\ar\quad\ar K\int_0^t\mu_{n-1}(t-s)J(\d s)\int_{\D_2}\int_{\D_2}\int_1^{\infty}\E f(\|\e^z{\bm R}(t,\mathbf{1})+{\bm u}+{\bm v}\|+1)\hat{m}_1(\d {\bm u})\hat{m}_2(\d{\bm v})\hat{\nu}(\d z).
\end{eqnarray*}
On the other hand,
\begin{eqnarray*}
\ar\ar\int_{\D_2}\int_{\D_2}\int_1^{\infty}\E f(\|\e^z{\bm R}(s,\mathbf{1})+{\bm u}+{\bm v}\|+1)\hat{m}_1(\d{\bm u})\hat{m}_2(\d{\bm v})\hat{\nu}(\d z)\cr\cr
\ar\le\ar  Kf(4)\int_{\D_2}\int_{\D_2}\int_1^{\infty}\E f(\frac{1}{4}\|\e^z{\bm R}(s,\mathbf{1})+{\bm u}+{\bm v}\|+1)\hat{m}_1(\d{\bm u})\hat{m}_2(\d {\bm v})\hat{\nu}(\d z)\cr\cr
\ar\le\ar \frac{1}{4}Kf(4)\bigg[K\int_1^{\infty}f(\e^z)\hat{\nu}(\d z)\E f\big(\|{\bm R}(s,\mathbf{1})\|\big)+\int_{\D_2}f(\|\z\|)\hat{m}_1(\d\z)\cr\cr
\ar\ar\qquad\qquad\qquad\qquad\qquad\qquad\qquad\qquad\qquad+\int_{\D_2}f(\|\z\|)\hat{m}_2(\d \z)+f(1)\bigg]\cr\cr
\ar\le\ar \frac{1}{4}Kf(4)\bigg[\sup\limits_{0\le t\le T}\E f(\|{\bm R}(t,\mathbf{1})\|)K\int_1^{\infty}f(\e^z)\hat{\nu}(\d z)\cr
\ar\ar\qquad\qquad\qquad\qquad\qquad\qquad\qquad\qquad+\int_{\D_2}f(\|\z\|)(\hat{m}_1+\hat{m}_2)(\d \z)+f(1)\bigg].
\end{eqnarray*}
Let $O_1(T)=\sup\limits_{0\le t\le T}\E f(\|{\bm R}(t,\mathbf{1})\|),$ $$O_2(T):=\frac{1}{4}K^2f(4)\Big[O_1(T)K\int_1^{\infty}f(\e^z)\hat{\nu}(\d z)+\int_{\D_2}f(\|\z\|)(\hat{m}_1+\hat{m}_2)(\d\z)+f(1)\Big].$$
The finiteness of $O_1(T)$ and $O_2(T)$ follow from Lemma \ref{lem2020020401}, $\int_{\{\|\z\|\ge 1\}}f(\|\z\|)(m_1+m_2)(\d\z)<\infty$ and $\int_1^{\infty}f(\e^z)\nu(\d z)<\infty.$ Then we have for any $t\in[0,T],$
$$
  \mu_n(t)\le O_1(T)+O_2(T)\int_0^t\mu_{n-1}(t-u)J_{\mathbf{1}}(\d u).
$$
On the other hand, according to Lemma 2 on pp.145 of \cite{ath}, there exists a positive function $t\mapsto \mu(t)$ bounded on $[0,T]$ s.t.
$$
  \mu(t)= O_1(T)+O_2(T)\int_0^t\mu(t-u)J_{\mathbf{1}}(\d u).
$$
Recall that $\mu_0(t)=0.$ It is clear that $\mu_n(t)\le\mu(t)$ for all $n=0,1,2,\cdots$ and $0\le t\le T.$ As $n$ tends to infinity, we have
$$\E f(\|\X(t,\x)\|)\le\lim\limits_{n\rightarrow\infty}\mu_n(t)\le\mu(t).$$ Using Proposition \ref{prop2020012901}, it is obvious that for any $\x\in\mathbb{R}_+^2,$ $t\mapsto\E f(\|\X(t,\x)\|)$ is locally bounded on $[0,\infty).$
\qed

%

\bigskip\bigskip

{\it Proof of Theorem 1.2:}

The sufficiency comes from Proposition \ref{prop2020012902} and  Proposition \ref{prop2020020402}. Conversely, if for some $t>0,$ $\E f(\|\X(t)\|)<\infty.$ Let $J_1(\d t)=\E (\rho_1\in\d t)$, where
$$\rho_1 = \inf \left\{t > 0: [X_1(t) - X_1(t-)]\wedge [X_2(t) - X_2(t-)]> 1, \xi_t = \xi_{t-}\right\}.
$$
Clearly, $\E\{\mathbf{1}_{\{\rho_1<t\}}\E[f(\|\X(t)\|)\vert\mathcal{F}_{\rho_1}]\}\leq\E f(\|\X(t)\|).$ By the strong Markov property,
 \begin{eqnarray*}
\ar\ar\E\{\mathbf{1}_{\{\rho_1<t\}}\E[f(\|\X(t)\|)\vert\mathcal{F}_{\rho_1}]\}\cr
\ar=\ar\E\{\mathbf{1}_{\{\rho_1<t\}}\mbb{E}_{\X(\rho_1)}[f(\|\x(t-\rho_1)\|)]\}\cr
\ar\ge\ar\E\{\mathbf{1}_{\{\rho_1<t\}}\mbb{E}_{\Delta\X(\rho_1)}[f(\|\x(t-\rho_1)\|)]\}\cr
\ar\ge\ar\E\int_0^t\{\mathbf{1}_{\{\rho_1<t\}}J_1(\d s)\int_{\mathbb{R}_+^2\setminus\{\mathbf{0}\}}\mbb{E}_{\z}f(\|\x(t-s)\|)\hat{m}_1(\d\z)\}\cr
\ar=\ar\E\{\int_0^t\mathbf{1}_{\{\rho_1<t\}}J_1(\d s)\int_{\mathbb{R}_+^2\setminus\{\mathbf{0}\}}\E f(\|\X(t-s,\z)\|)\hat{m}_1(\d\z)\}\cr
\ar=\ar\E\{\int_0^t\mathbf{1}_{\{\rho_1<t\}}J_1(\d s)\int_{\mathbb{R}_+^2\setminus\{\mathbf{0}\}}\E f(\|{\bm R}(t-s,\z)\|)\hat{m}_1(\d\z)\},
\end{eqnarray*}
In the above equations, $\{{\bm R}(t,\z):t\ge 0\}$ is the strong solution of \eqref{2020020406}-\eqref{2020020407} with initial value $\z.$ Since $t\mapsto J_1(0,t]$ is strictly increasing, there exists some $s\in(0,t],$ such that
$$
\int_{\mathbb{R}_+^2\setminus\{\mathbf{0}\}}\E f(\|{\bm R}(t-s,\z)\|)\hat{m}_1(\d\z)<\infty.
$$
By Proposition \ref{prop2020012902}, $\int_{\D_2} f(\|\z\|)\hat{m}_1(\d\z)<\infty.$
Then $\int_{\D_2}f(\|\z\|)m_1(\d\z)<\infty.$ Similarly, $\int_{\D_2}f(\|\z\|)m_2(\d\z)<\infty.$
In conclusion, $\int_{\{\|\z\|\ge 1\}}f(\|\z\|)(m_1+m_2)(\d\z)<\infty.$ Similarly, define  $\rho_2 = \inf \left\{t > 0: \xi(t) - \xi(t-) > 1\right\},$ $J_2(\d t)=\p (\rho_2\in\d t).$ By the strong Markov property,
\begin{eqnarray*}
\ar\ar\E\{\mathbf{1}_{\{\rho_2<t\}}\E[f(\|\X(t)\|)\vert\mathcal{F}_{\rho_2}]\}\cr
\ar=\ar\E\{\mathbf{1}_{\{\rho_2<t\}}\mbb{E}_{\X(\rho_2)}[f(\|\x(t-\rho_2)\|)]\}\cr
\ar\ge\ar\int_1^{\infty}\E\{\mathbf{1}_{\{\rho_2<t\}}\mbb{E}_{\e^z {\bm R}(\rho_2)}[f(\|\x(t-\rho_2)\|)]\}\hat{\nu}(\d z)\cr
\ar\ge\ar\int_0^t\mathbf{1}_{\{{\bm R}(\rho_2)>{\bm \varepsilon}\}}J_2(\d s)\int_1^{\infty}\mbb{E}_{\e^z{\bm\varepsilon}}f(\|\x(t-s)\|)\hat{\nu}(\d z)\cr
\ar=\ar\int_0^t\mathbf{1}_{\{{\bm R}(\rho_2)>{\bm \varepsilon}\}}J_2(\d s)\int_1^{\infty}\E f(\|\X(t-s,\e^z{\bm\varepsilon})\|)\hat{\nu}(\d z)\cr
\ar\ge\ar\int_0^t\mathbf{1}_{\{{\bm R}(\rho_2)>{\bm\varepsilon}\}}J_2(\d s)\int_1^{\infty}\E f(\|{\bm R}(t-s,\e^z{\bm\varepsilon})\|)\hat{\nu}(\d z),
\end{eqnarray*}
 where $\{{\bm R}(t,\e^z{\bm\varepsilon}):t\ge 0\}$ is the strong solution of \eqref{2020020406}--\eqref{2020020407} starting from $\e^z{\bm\varepsilon}$ with ${\bf\varepsilon}\in\mathbb{R}^2_+\setminus\{\mathbf{0}\}$ s.t. $\p({\bm R}(t-s)>{\bm\varepsilon})>0.$ Thus $t\mapsto J_2(0,t]$ is strictly increasing. Therefore, there exists $s\in(0,t]$ such that
 $$
\int_1^{\infty}\E f\Big(\|{\bm R}(t-s,\e^z{\bm\varepsilon})\|\Big)\hat{\nu}(\d z)<\infty.
$$
By Proposition \ref{prop2020012902},
$$\int_1^{\infty}f(\e^z)\hat{\nu}(\d z)\le Kf(\|{\bm\varepsilon}\|^{-1})\int_1^{\infty}f(\|\e^z{\bm\varepsilon}\|)\hat{\nu}(\d z).$$
Hence, $\int_1^{\infty}f(\e^z)\nu(\d z)<\infty$ and $\E f(\|\X(0)\|)<\infty.$
\qed
\section{Future research}

In this paper we calculate the integer moments and $f$-moment for fixed $t\ge 0,$ where $f$ is a function on $[0, \infty).$ The case when $f$ is a bivariate function and the moment behavior as $t\rightarrow\infty$ is still to be explored. Also, some probability distributions can be uniquely determined by the integeter moments, see \cite{st1} and \cite{st2} for instance. And we are interested in the moment determinacy of CB-processes for fixed $t\ge 0$ and also the limit as $t\rightarrow\infty.$ Furthermore, since we are discussing two-type CBRE-process here, another direction for the future is to allow the number of types going up to countable infinity, like   the authors in \cite{kyp} did.
\section*{Acknowledgments}

Heartfelt thanks are given to Professor Jordan Stoyanov for his valuable comments on this work. The authors also thank Doctor Lina Ji and Doctor Rongjuan Fang for their useful discussions and suggestions. This paper is supported by Special Fund for Central Universities SLK13223001.

\section*{Declarations}
\begin{itemize}
\item The authors have no conflicts of interest to declare. All co-authors have seen and agree with the contents of the manuscript.
\item  The data that support the findings of this study are available on request from the corresponding author.
\end{itemize}


\begin{thebibliography}{00}
\bibitem{ath}Athreya, K. B. and Ney, P. E. {\it Branching Processes.} Berlin: Springer, 1972.
\bibitem{bar}Barczy, M., Li, Z. and Pap, G. Stochastic differential equation with jumps for multi-type continuous state and continuous time branching processes with immigration. {\it ALEA Lat Am J Probab Math Stat.}, 2015, \textbf{12}(1):129-169.
\bibitem{barm}Barczy, M., Li, Z. and Pap, G. Moment formulas for multitype continuous state and continuous time branching process with immigration. {\it J. Theor. Probab.}, 2016, \textbf{29}(3): 958-995.
\bibitem{ref8}
He, H., Li, Z. and Xu, W. Continuous-State Branching Processes in L\'{e}vy Random Environments. {\it J. Theor. Probab.}, 2018, \textbf{31}: 1952-1974.
\bibitem{Ji58}Ji\u{r}ina M. Stochastic branching processes with continuous state space. {\it Czechoslov. Math. J.}, 1958, \textbf{8}: 292-313.
\bibitem{ref6}
Ji, L. and Li, Z. Moments of Continuous-state Branching Processes with or Without Immigration. {\it Acta Math Appl Sin Engl Ser.}, 2020, \textbf{36}, 361-373.
\bibitem{ref7}
Ji, L. and Zheng, X. Moments of Continuous-State Branching Processes in L\'{e}vy Random Environments. {\it Acta Math. Sci. Engl. Ser.,} 2019, \textbf{39}: 781-796.
\bibitem{ref2}
Kawazu, K. and Watanabe, S. Branching processes with immigration and related limit theorems. {\it Theory Probab. its Appl.}, 1971, \textbf{16}: 36-54.
\bibitem{kyp}
Kyprianou, A. E. and Palau S. Extinction properties of multi-type continuous-state branching processes. {\it Stoch Process Their Appl}, 2018, \textbf{128}(10): 3466-3489.
\bibitem{ref3}
Li, Z. A limit theorem for discrete Galton-Watson branching processes with immigration. {\it J Appl Probab}, 2006, \textbf{43}: 1103-1142.
\bibitem{ref4}
Li, Z. {\it Measure-Valued Branching Markov Processes.} Heidelberg: Springer, 2011.
\bibitem{machunhua}Ma, C. A limit theorem of two-type Galton-Watson branching processes with immigration.{\it Stat Probab Lett}, 2009, \textbf{79}, no. 15, 1710-1716.
\bibitem{marugang}Ma, R. Stochastic equations for two-type continuous-state branching processes
with immigration and competition. {\it Stat Probab Lett}, 2014, \textbf{91}: 83-89.
\bibitem{ref9}
Palau, S. and Pardo, J C. Branching processes in a L\'{e}vy random environment. {\it Acta Appl Math}, 2018, \textbf{153}(1): 55-79.
\bibitem{qinzheng}
Qin, Y. and Zheng, X. Stochastic equations and ergodicity for two-type continuous-state branching processes with immigration in L\'evy random environments. {\it Math. Methods Appl. Sci.}, 2020, \textbf{43}: 8363-8378.
\bibitem{Sato}
Sato, K I. {\it L\'{e}vy Processes and Infinitely Divisible Distributions}. Cambridge university press, 1999.
\bibitem{st1}Stoyanov, J. M, Lin, G. D. and Kopanov, P. New checkable conditions for moment determinacy of probability distributions. {\it Theory Probab. its Appl.} 2020, 65(3): 497-509.
\bibitem{st2}Stoyanov J. M. Moment properties of probability distributions used in stochastic financial models. RECENT ADVANCES IN FINANCIAL ENGINEERING 2014: Proceedings of the TMU Finance Workshop 2014. 2016: 1-27.
\end{thebibliography}
\end{document}